\documentclass[10pt]{article}
\usepackage{graphicx,amsmath,amssymb,a4wide}

\begin{document}
\newtheorem{thme}{Theorem}[section]
\newtheorem{prop}[thme]{Proposition}
\newtheorem{lemme}[thme]{Lemma}
\newtheorem{rem}[thme]{Remark}
\newtheorem{defin}[thme]{Definition}
\newtheorem{cor}[thme]{Corollary}
\newtheorem{question}[thme]{Question}
\newcommand{\R}{\mathbb{R}}
\newcommand{\A}{\mathbb{A}}
\newcommand{\B}{\mathbb{B}}
\newcommand{\C}{\mathbb{S}^1}
\newcommand{\At}{\widetilde{\mathbb{A}}}
\newcommand{\Ab}{\breve{\mathbb{A}}}
\newcommand{\Bb}{\breve{\mathbb{B}}}
\newcommand{\N}{\mathbb{N}}
\newcommand{\Z}{\mathbb{Z}}
\newcommand{\cqfd}{\hfill $\blacksquare$}

\title{A topological version of the Poincar\'e-Birkhoff theorem with two fixed points}
\author{Marc Bonino}
\date{}
\maketitle 
\begin{center}
\small  
Laboratoire Analyse, G\'eom\'etrie et Applications (LAGA) \\
CNRS UMR 7539 \\ 
Universit\'e Paris 13 \\ 
99 Avenue J.B. Cl\'ement \\
93430 Villetaneuse (France) \\
e-mail: \textup{\texttt{bonino@math.univ-paris13.fr}} 
\end{center}
\textbf{Abstract.} The main result of this paper gives a topological property 
satisfied by any homeomorphism of the annulus $\A=\C \times [-1,1]$ isotopic to the identity and with at most one fixed point. This generalizes the classical Poincar\'e-Birkhoff theorem because this property certainly does not hold for an area preserving homeomorphism $h$ of $\A$ with the usual boundary twist condition. We also have two corollaries of this result. The first one shows in particular that the boundary twist assumption may be weakened by demanding that the 
homeomorphism $h$ has a lift $H$ to the strip $\At = \R \times [-1,1]$ possessing both a forward orbit unbounded on 
the right and a forward orbit unbounded on the left. As a second corollary we get a new proof of a version of the Conley-Zehnder theorem in $\A$: if a homeomorphism of $\A$ isotopic to the identity preserves the area and has mean rotation zero, then it possesses two fixed points.  

\normalsize  
\section{Introduction} \label{section1}
In its more classical version, the Poincar\'e-Birkhoff theorem is the following.
\begin{thme}[Birkhoff, \cite{Birkhoff:1913,Birkhoff:1925}] \label{tPB} 
Let $h$ be a homeomorphism of the annulus $\A = \mathbb{S}^1 \times [-1, 1]$ isotopic to the identity.
If $h$ preserves the area and satisfies the boundary twist condition then it has at least two fixed points.
\end{thme}
The boundary twist condition requires that $h$ has a lift $H$ to the universal cover 
$\widetilde{\A} = \R \times [-1, 1]$ such that, writing $H(\theta, \pm 1) = (\varphi_\pm (\theta), \pm 1)$, one has $(\varphi_- (\theta) - \theta)(\varphi_+ (\theta)-\theta) < 0$ for every $\theta \in \R$. 
Let us recall that this result was conjectured by Poincar\'e in \cite{Poincare:1912} 
and proved by Birkhoff in \cite{Birkhoff:1913,Birkhoff:1925} (see also \cite{Brown/Neumann:1977}). In fact the proof in \cite{Birkhoff:1913} actually ensures the existence of only one fixed point while the full Theorem \ref{tPB} turns out to be a particular case of a more general statement in \cite{Birkhoff:1925} where the annulus is not necessarily setwise invariant under the homeomorphism. This is an interesting direction for generalizing Theorem \ref{tPB} (in particular for applications to ODE's, see e.g. \cite{Dalbono/Rebelo:2002}) but we focus in this paper on self-homeomorphisms of $\A$. Our interest is in the search of more topological versions and on the number of fixed points that one can expect from such variants. These questions seem to originate from \cite{Poincare:1912}; indeed it is well known that Poincar\'e overlooked the possibility for a homeomorphism of $\A$ to have a single fixed point and so suggested the following strategy for proving his theorem (\cite{Poincare:1912}[p.376-377]): assume that $T$ is a homeomorphism of $\A$ satisfying the boundary twist condition and without fixed point; then one could show that $T$ does not preserve the area by constructing an essential Jordan curve $C \subset \A$ disjoint from its image $C'=T(C)$. Precisely this program was achieved by Ker\'ekj\'art\'o in \cite{Kerekjarto:1928} and the right generalization for obtaining the second fixed point was obtained by P. H. Carter in \cite{Carter:1982} \footnote{As Birkhoff in \cite{Birkhoff:1925}, Carter deals with the situation where the annulus is not setwise invariant.}. Other references for Ker\'ekj\'art\'o's and Carter's results are \cite{Franks:1988b}, \cite{Guillou:1997} and \cite{LeCalvez/Wang:2009}. One can loosely say after Carter's work that the area preserving hypothesis in Theorem \ref{tPB} may be replaced with the weaker assumption that there is no essential subannulus $\B \subset \A$ containing its image as a proper subset. 
Afterwards some authors generalized both the twist and the conservative hypotheses; one can quote C. Bonatti and L. Guillou (\cite{Guillou:1994}[Th\'eor\`eme 5.1]), J. Franks (\cite{Franks:1988a,Franks:1988c}) and H.E. Winkelnkemper (\cite{Winkelnkemper:1988}), the results in \cite{Franks:1988c}, \cite{Guillou:1994} and \cite{Winkelnkemper:1988} giving only one fixed point. 

The first result of the present paper is a purely topological version of the Poincar\'e-Birkhoff theorem allowing to detect two fixed points (Theorem \ref{t1}). 
More precisely this result shows that a homeomorphism $h:\A \to \A$ isotopic to identity and possessing at most one fixed point should be either ``non-conservative" or ``untwisted" in some topological sense. It appears as a natural extension of 
Bonatti-Guillou's theorem concerning fixed point free homeomorphism of $\A$. Section \ref{section5} then relates Theorem \ref{t1} (and its refinement Theorem \ref{t2}) with other works on the subject. As a first application, we obtain Theorem \ref{t3} showing that a homeomorphism $h:\A \to \A$ can be thought of as twisted if there 
is a lift $H:\At \to \At$ possessing a forward orbit unbounded on the left and 
another one unbounded on the right. This extends a theorem by Franks where the twist 
assumption is interpreted by means of rotation numbers (\cite{Franks:1988a}[Theorem 3.3]). We get finally a short ``geometric" proof of M. Flucher's version of the Conley-Zehnder theorem in the 
annulus (\cite{Flucher:1990}[Theorem 2]) when dealing with homeomorphism preserving the Lebesgue 
measure (Theorem \ref{t4}). 

\ \\
This article can also be regarded as the continuation of \cite{Bonino:2009} where results close to the 
Poincar\'e-Birkhoff theorem were obtained in the isotopy class of the symmetry $S_{\A}$ interchanging the boundary components of $\A$. Theorem \ref{t2} below and Theorem 1.2 of \cite{Bonino:2009} have  
indeed strong similarities, as well as their proofs. Nevertheless these two results are also logically independent 
and the present paper is largely self-contained.

\section{Preliminaries} \label{section2}
\subsection{Definitions and notation} 
The notation and vocabulary used throughout this paper are the same as in \cite{Bonino:2009}. In particular 
the strip $\At = \R \times [-1,1]$ is regarded as the universal cover of the annulus $\A= \C \times [-1,1]$ with covering map $\Pi(\theta,r)= (e^{2i\pi\theta},r)$. The deck transformations are then the iterates $\tau^n: \At \to \At$ 
of the translation $\tau(\theta,r)=(\theta+1,r)$. 
We write $\overline{X}, \mathrm{Int}(X), \partial X$ for respectively the closure, interior and frontier relative to 
$\At$ of a subset $X \subset \At$. When we need to consider these notions with respect to another topological space $Y$, 
we use the explicit notation $\mathrm{Cl}_Y(X), \mathrm{Int}_Y(X)$ and  $\partial_Y X$ for any $X \subset Y$. 
The reader is refered to the short Section 2.1 of \cite{Bonino:2009} for other details.
Recall also that if $f:\R^2 \to \R^2$ is a continuous map and $J$ a Jordan curve in $\R^2 \setminus \mathrm{Fix}(f)$, 
then the \emph{index} of $J$ w.r.t. $f$ is the winding number of the vector $f(p)-p$ when $p$ describes 
$J$ counterclokwised oriented.  Given any map $f:E \to E$, a family $\mathcal{E}$ of subsets of $E$ is said to 
be $f$-\emph{free} if $f(X) \cap X = \emptyset$ for every $X \in \mathcal{E}$. Finally a subset $\B \subset \A$ is named a \emph{subannulus} of $\A$ if it is homeomorphic to $\A$. 

\subsection{Nielsen classes} If $h:\A \to \A$ is a continuous map, its Nielsen classes may be defined as the various sets $\Pi(\mathrm{Fix}(H))$ where $H:\widetilde{\A} \to \widetilde{\A}$ varies over all the lifts of $h$ (see e.g. \cite{Jiang:1983}). 
The nonempty Nielsen classes of $h$ realize a finite partition of $\mathrm{Fix}(h)$. 
We consider throughout this paper only maps $h:\A \to \A$ which are homotopic to the identity; the Nielsen classes $\Pi(\mathrm{Fix}(H))$ and $\Pi(\mathrm{Fix}(H'))$ are then disjoint whenever $H,H': \At \to \At$ are two distinct lifts of $h$: this follows easily from the fact that one has in this setting $H \circ \tau = \tau \circ H$ and since $H' = \tau^n \circ H$ for some $n \in \Z \setminus \{0\}$. We will say that two Nielsen classes $N,N'$ are \emph{consecutive} if they are both nonempty and if, for some $k \in \{\pm 1\}$ and for some lift $H: \At \to \At$ of $h$, one has $N = \Pi(\mathrm{Fix}(H))$ and $N' =\Pi(\mathrm{Fix}(\tau^k \circ H))$. 

\subsection{Le Calvez-Sauzet's brick decompositions} 
This notion was introduced in \cite{LeCalvez/Sauzet:1996,Sauzet:2001} as a convenient tool for proving Brouwer's plane translation theorem. A \emph{brick decomposition} of a surface $S$ consists essentially in a locally finite tiling of $S$ with topological closed discs (the \emph{bricks} of the decomposition) in such a way that any subset $X \subset S$ 
obtained as the union of some bricks is a subsurface of $S$. 
We refer to Section 2.4 of \cite{Bonino:2009} for a definition and for basic properties. 
Because of their importance, we recall now the related notions of attractor and repeller. 
Suppose that $f: S \to S$ is a homeomorphism of a surface $S$ endowed with a brick decomposition 
$\mathcal{D} = \{B_i\}_{i \in I}$. The \emph{attractor} associated to a brick $B_{i_0}$ (and to $f$) is the union of 
all the bricks $B$ such that, for some integer $n \geq 1$, there exists a sequence of bricks $B_{i_0}, B_{i_1}, \cdots,B_{i_n}=B$ 
satisfying $$\forall j \in \{0,\cdots,n-1\} \, f(B_{i_j}) \cap B_{i_{j+1}} \neq \emptyset.$$ One defines 
similarly the $\emph{repeller}$ associated to $B_{i_0}$ by replacing $f$ with $f^{-1}$.  
If $\mathcal{A,R}$ denote respectively the attractor and the repeller associated to some brick, one 
clearly has $f(\mathcal{A}) \subset \mathcal{A}$ and $f^{-1}(\mathcal{R}) \subset \mathcal{R}$. Even better, since the 
union of the bricks containing a given point $z \in S$ is a neighbourhood of $z$ in $S$, one has $f(\mathcal{A}) \subset \mathrm{Int}_S(\mathcal{A})$ and $f^{-1}(\mathcal{R}) \subset \mathrm{Int}_S(\mathcal{R})$. We also have 
the following essential properties.
 \begin{prop}\label{p1}
Let $f:\R^2 \to \R^2$ be an orientation preserving homeomorphism leaving setwise invariant 
a subsurface $S \subset \R^2 \setminus \mathrm{Fix}(f)$ and such that no Jordan curve 
$J \subset \R^2 \setminus \mathrm{Fix}(f)$ has index $1$ w.r.t. $f$. Let $B_{i_0}$ be any brick of a $f$-free 
brick decomposition $\mathcal{D}$ of $S$ and write $\mathcal{A,R}$ for respectively the attractor and the repeller associated to $B_{i_0}$ and $f$. Then  
\begin{itemize}
\item [(i)] $B_{i_0}$ is contained neither in $\mathcal{A}$ nor in $\mathcal{R}$, i.e. 
$B_{i_0} \cap  \mathrm{Int}_S(\mathcal{A}) = \emptyset = B_{i_0} \cap \mathrm{Int}_S(\mathcal{R})$;
\item[(ii)] there is no brick contained in $\mathcal{A} \cap \mathcal{R}$, i.e. $\mathrm{Int}_S(\mathcal{A}) 
\cap \mathcal{R} = \emptyset = \mathrm{Int}_S(\mathcal{R}) \cap \mathcal{A}$.
\end{itemize}
\end{prop}
This follows from Franks' lemma about periodic disc chains
(\cite{Franks:1988a}[Proposition 1.3]) which is itself a consequence of Brouwer's lemma on translation arcs. Nevertheless one should observe that Franks' statement deals with open (topological) discs whereas we use it with closed discs, namely the bricks of $\mathcal{D}$. The needed refinements to work with closed discs are due to Guillou and Le Roux (see \cite{LeRoux:2004}[p.38-39]).

\section{Statement of the main results} \label{section3}
Our first result is the following.
\begin{thme} \label{t1}
Let $h:\A \to \A$ be a homeomorphism of the compact annulus $\A = \mathbb{S}^1 \times [-1,1]$ isotopic to the identity 
(i.e. $h$ preserves the orientation and the two boundary components of $\A$). If $h$ has no more than one fixed point then at least one of the two following properties holds:
\begin{enumerate}
\item There exists an essential Jordan curve $J \subset \A$ such that 
$J \cap h(J) = J \cap \mathrm{Fix}(h)$. 
\item There exists an arc $\alpha$ crossing $\A$ such that
\begin{itemize}
\item $\alpha \cap h(\alpha) = \alpha \cap \mathrm{Fix}(h)$,
\item $h(\alpha)$ does not meet the two local sides of $\alpha$.
\end{itemize}
\end{enumerate}
\end{thme}
As mentioned in the introduction, Theorem \ref{t1} is a natural extension of Bonatti-Guillou's 
version of the Poincar\'e-Birkhoff theorem since it reduces exactly to \cite{Guillou:1994}[Th\'eor\`eme 5.1] 
under the strongest assumption that $h$ is fixed point free.  According to the Lefschetz-Hopf theorem, 
if $z \in \A$ is the unique fixed point of a continuous map $h:\A \to \A$ homotopic to identity then its Lefschetz index 
equals the Euler characteristic $\chi(\A)=0$ hence Theorem \ref{t1} is contained in Theorem \ref{t2} below. 
This result points out, we hope, what are the truly important assumptions for our constructions and will be useful for 
Theorem \ref{t3} and Theorem \ref{t4}. 
\begin{thme} \label{t2}
Let $h:\A \to \A$ be a homeomorphism of the compact annulus $\A = \mathbb{S}^1 \times [-1,1]$ isotopic 
to the identity. If $h$ satisfies the following three assumptions
\begin{enumerate} 
\item[(i)] There is at least one boundary component of $\A$ where $h$ has no fixed point;
\item[(ii)] There are no two consecutive Nielsen classes;  
\item[(iii)] $h$ only has finitely many fixed points and these fixed points (if any) have Lefschetz index $0$; 
\end{enumerate}
then, for some Nielsen class $N_0$ of $h$ (maybe $N_0=\emptyset$), at least one of the following two properties holds:
\begin{enumerate}
\item[1'.] There exists an essential subannulus $\B \subset \A$ containing either $h(\B)$ or 
$h^{-1}(\B)$ as a proper subset. Moreover the boundary of $\B$ consists of one of the 
two boundary component of $\A$, call it $\mathrm{Bd}^{\sigma}(\A)$, together with an essential Jordan curve 
$J \subset \A \setminus \mathrm{Bd}^{\sigma}(\A)$ such that $J \cap h(J) = J \cap N_0$;
\item[2'.] There exists an arc $\alpha$ crossing $\A$ such that
\begin{itemize}
\item $\alpha \cap h(\alpha) = \alpha \cap N_0$,
\item $h(\alpha)$ does not meet the two local sides of $\alpha$.
\end{itemize}
\end{enumerate}
\end{thme}

\section{Proof of Theorem \ref{t2}} \label{section4}
Let $h: \A \to \A$ be a homeomorphism isotopic to the identity and satisfying the assumptions (i)-(iii) in Theorem \ref{t2}. 
We can suppose without loss that $h$ is fixed point free 
on $\mathrm{Bd}^-(\A)$ so there exists a lift $H_0$ of $h$ to $\At$ such that 
$$ \forall (\theta,-1) \in \mathrm{Bd}^{-}(\At) \quad (\theta,-1) < H_0(\theta,-1) < (\theta + 1,-1) $$ 
and then also 
$$ \forall (\theta,-1) \in \mathrm{Bd}^{-}(\At) \quad (\theta,-1) < \tau \circ H_0^{-1}(\theta,-1) < (\theta+1,-1).$$
According to $(ii)$ we can pick $G \in \{H_0,\tau \circ H_0^{-1} \}$ which is fixed point free and we rename $H$ the remaining homeomorphism. Thus $H$ is a lift of $h$ while $G$ is a lift of $h^{-1}$ or conversely, and anyway 
$G \circ H = H \circ G = \tau$. We let $F=\mathrm{Fix}(H)$ and $S= \At \setminus F$. Thus $S$ is a subsurface of $\At$ such that $\mathrm{Bd}^-(\At) \subset S = \tau(S)=H(S)$ and it is an open subset of $\At$. We also define $N_0$ to be the Nielsen class of $h$ determined by $H$, that means $N_0=\Pi(F)$. We have the following easy fact which will allows us to use Proposition \ref{p1}.
\begin{lemme} \label{l0}
The homeomorphism $H:\At \to \At$ can be extended to an orientation preserving homeomorphism $f:\R^2 \to \R^2$ such 
that any Jordan curve $J \subset \R^2 \setminus \mathrm{Fix}(f)$ have index $0$ w.r.t. $f$. 
\end{lemme}
\emph{Proof of Lemma \ref{l0}.} Just let, for $(x,y) \in \R^2$ and $\vert y \vert \geq 1$,
$$f(x,y) = H(x,\dfrac{y}{\vert y \vert}) + (0,y- \dfrac{y}{\vert y \vert})$$
and remark that any connected component of $\mathrm{Fix}(f)$ is either (vertically) unbounded or consists in a fixed point of $H$ in $\At \setminus \mathrm{Bd}(\At)$. Due to $(iii)$, all the fixed points of $H$ are isolated with index $0$ and 
the result follows since the index of a Jordan curve $J \subset \R^2 \setminus \mathrm{Fix}(f)$ surrounding finitely 
many fixed points $p_1, \cdots, p_n$ is equal to the sum of the indices of the $p_i$'s. \cqfd 

\ \\
We construct now a brick decomposition of $S$ adapted to our purpose. In particular the fourth item in Lemma \ref{l1}
ensures conveniently that the ``dynamics on the bricks" adjacent to $\mathrm{Bd}^-(\At)$ looks like the behaviour of $H$ on 
$\mathrm{Bd}^-(\At)$. Similar ideas where used for Lemma 3.1 of \cite{Bonino:2009}. 
\begin{lemme} \label{l1} There exist an $\epsilon \in (0,1)$ and a brick decomposition ${\widetilde{\mathcal{D}}}_H = \{B_i\}_{i \in I}$ of the surface~$S$ satisfying the following properties:
\begin{enumerate}
\item $\widetilde{\mathcal{D}}_H$ is $\tau$-equivariant, which means $\{\tau(B_i)\}_{i \in I} = {\widetilde{\mathcal{D}}}_H$.
\item ${\widetilde{\mathcal{D}}}_H$ is $H$-free.
\item The bricks meeting $\mathrm{Bd}^-(\At)$ are rectangles $B^-_n = [a_n,a_{n+1}] \times [-1,-1+\epsilon]$ where $(a_n)_{n \in \mathbb{Z}}$ is a strictly 
increasing sequences of reals numbers such that $\lim_{\pm \infty} a_n = \pm \infty$. 
\item For every $n \in \mathbb{Z}$ one has $H(B^-_n) \cap B^-_{n+1} \neq \emptyset$, hence $B^-_{n+1} \subset \mathcal{A}(B^-_n)$ and $B^-_{n-1} \subset \mathcal{R}(B^-_n)$.
\end{enumerate}
\end{lemme}
\emph{Proof of Lemma \ref{l1}.} A family $\mathcal{I}$ of compact intervals of $\mathrm{Bd}^-(\A)$ is named an \emph{interval decomposition} of $\mathrm{Bd}^-(\A)$ if $\bigcup_{\alpha \in \mathcal{I}} \alpha 
= \mathrm{Bd}^-(\A)$ and if any two distinct $\alpha,\alpha'$ of $\mathcal{I}$ can meet only in a common endpoint. We define similarly an interval decomposition of $\mathrm{Bd}^-(\At)$. Because $h$ has no fixed point 
on $\mathrm{Bd}^-(\A)$ there exists a finite and $h$-free interval decomposition $\mathcal{I}_0$ of $\mathrm{Bd}^-(\A)$ (with at least three intervals). By considering all the connected components of all the sets $\Pi^{-1}(\alpha)$, $\alpha \in \mathcal{I}_0$, one gets a $H$-free and $\tau$-equivariant interval decomposition $\widetilde{\mathcal{I}_0}$ of $\mathrm{Bd}^-(\At)$. Given two interval decompositions $\widetilde{\mathcal{I}}, \widetilde{\mathcal{I}}'$ of $\mathrm{Bd}^-(\At)$, let us 
write $\widetilde{\mathcal{I}}' \succeq \widetilde{\mathcal{I}}$ if every interval of $\widetilde{\mathcal{I}}$ is contained in an interval of $\widetilde{\mathcal{I}}'$. This defines a partial ordering on the set of all the interval decompositions of $\mathrm{Bd}^-(\At)$ and there exists $\widetilde{\mathcal{I}} \succeq \widetilde{\mathcal{I}_0}$ which is maximal among the interval decompositions of $\mathrm{Bd}^-(\At)$ which are both $H$-free and $\tau$-equivariant. 
Observe that if $\tilde{\alpha} = [a,a'] \times \{-1\}$ and $\tilde{\alpha}' = [a',a''] \times \{-1\}$ are two consecutive intervals of $\widetilde{\mathcal{I}}$ ($a<a'<a''$) then necessarily $H(\tilde{\alpha}) \cap \tilde{\alpha}' \neq \emptyset$. Indeed one has $\tau(\tilde{\alpha}) \neq \tilde{\alpha}'$ because 
$H(\tilde{\alpha}) \cap \tilde{\alpha} = \emptyset$ and $(a,-1)< H(a,-1) < (a+1,-1)$ on $\mathrm{Bd}^-(\At)$ so one gets another $\tau$-equivariant interval decomposition $\widetilde{\mathcal{I}}' \succeq \widetilde{\mathcal{I}}$ by 
removing from $\widetilde{\mathcal{I}}$ all the intervals $\tau^n(\tilde{\alpha}), \tau^n(\tilde{\alpha}')$ and 
replacing them with the $\tau^n( \tilde{\alpha} \cup \tilde{\alpha}')$'s ($n \in \Z$). Due to the maximality of 
$\widetilde{\mathcal{I}}$, this $\widetilde{\mathcal{I}}'$ cannot be $H$-free so $H(\tilde{\alpha} \cup \tilde{\alpha}') \cap (\tilde{\alpha} \cup \tilde{\alpha}') \neq \emptyset$. Since $H$ moves the points of $\mathrm{Bd}^-(\At)$ towards the right one gets $H(\tilde{\alpha}) \cap \tilde{\alpha}' \neq \emptyset$, as announced. 
The family $\{\Pi(\tilde{\alpha})\}_{\tilde{\alpha} \in \widetilde{\mathcal{I}}}$ being finite, there exists $\epsilon \in (0,1)$ so small that, for every $\tilde{\alpha} =[a,a'] \times \{-1\} \in \widetilde{\mathcal{I}}$, the rectangle $\tilde{R}_{\tilde{\alpha}} = [a,a'] \times [-1,1+\epsilon]$ is disjoint from its $H$-image. 
This provides the bricks $B_i^-$ adjacent to $\mathrm{Bd}^-(\At)$ as described in (3)-(4) and it remains to complete 
this collection of rectangles to get the required brick decomposition of $S$. 
The family $\big \{ \Pi(\tilde{R}_{\tilde{\alpha}}) \big \}_{\tilde{\alpha} \in \widetilde{\mathcal{I}}}$ 
consists in finitely many topological closed discs in $\A \setminus N_0$ and one easily constructs a brick decomposition $\mathcal{D}$ of 
$\A \setminus N_0$ containing the $\Pi(\tilde{R}_{\tilde{\alpha}})$'s as bricks. Subdividing if necessary some bricks of $\mathcal{D}$ others than the $\Pi(\tilde{R}_{\tilde{\alpha}})$'s, one can assume that all the connected components of the sets $\Pi^{-1}(B)$, $B \in \mathcal{D}$, are disjoint from their $H$-image; hence they provide a brick decomposition of $S$ 
satisfying (1)-(4).\cqfd 

\ \\
We fix from now on a brick decomposition ${\widetilde{\mathcal{D}}}_H$ of $S$ given by Lemma \ref{l1} and a brick $B^-_i$ adjacent to $\mathrm{Bd}^-(\At)$, say $B^-_0 = [a_0,a_1] \times [-1,-1+\epsilon]$. We write respectively $\mathcal{A,R}$ for the attractor and the repellor associated to $B^-_0$ and $H$. Property (4) in Lemma \ref{l1} ensures 
that $[a_1,+ \infty) \times \{-1\} \subset \bigcup_{i \geq 1} B_i^- \subset \mathcal{A}$ so $\mathcal{A}$ is unbounded on the right and moreover, using $B_0^- \cap B_1^- \neq \emptyset$, that $\mathcal{A}$ is connected. Similarly 
$\mathcal{R}$ is connected and unbounded on the left, with 
$(- \infty, a_0] \times \{-1\} \subset \bigcup_{i \leq -1} B_i^- \subset \mathcal{R}$.
Theorem \ref{t2} now follows from Propositions \ref{p2}-\ref{p3} below. 
\begin{prop} \label{p2}
The following implications holds: 
\begin{itemize}
\item $\mathcal{A}$ is unbounded on the left $\Rightarrow \mathcal{R} \cap \mathrm{Bd}^+(\At)
= \emptyset \Rightarrow$ the alternative (1') of Theorem \ref{t2} occurs. 
\item $\mathcal{R}$ is unbounded on the right $\Rightarrow \mathcal{A} \cap \mathrm{Bd}^+(\At)
= \emptyset \Rightarrow$ the alternative (1') of Theorem \ref{t2} occurs. 
\end{itemize}
\end{prop}
\begin{prop} \label{p3}
If $\mathcal{A}$ is bounded on the left and meets $\mathrm{Bd}^+(\At)$ then the conclusion of Theorem \ref{t2} holds. 
\end{prop}
\emph{Proof of Proposition \ref{p2}.} 
We just prove the first point, the second one being similar by reversing the roles of $\mathcal{A}$ and $\mathcal{R}$. 
First suppose that $\mathcal{R} \cap \mathrm{Bd}^+(\At) \neq \emptyset$. The set $\mathrm{Int}(\mathcal{R})$ is 
arcwise connected (as the interior of a connected union of bricks) so there 
exists an arc $\tilde{\gamma}$ crossing $\At$ and entirely contained in $\mathrm{Int}(\mathcal{R})$. 
Proposition \ref{p1} then gives $\mathcal{A} \cap \tilde{\gamma} \subset \mathcal{A} \cap \mathrm{Int}(\mathcal{R}) 
= \emptyset$. Since furthermore $\mathcal{A}$ is unbounded on the right and connected, 
it should be contained in the domain on the right of $\tilde{\gamma}$ and so it is bounded on the left. 

Assume now that $\mathcal{R} \cap \mathrm{Bd}^+(\At) = \emptyset$ and define 
$X_0 = \bigcup_{n \in \mathbb{Z}} \tau^n(\mathcal{R}) \subset \At \setminus \mathrm{Bd}^+(\At)$. 
The brick decomposition $\widetilde{\mathcal{D}}_H$ being $\tau$-equivariant, each $\tau^n(\mathcal{R})$ is a union 
of bricks of $\widetilde{\mathcal{D}}_H$ and so is $X_0$. One clearly has $\bigcup_{i \in \Z} B_i^- \subset X_0$.  
Moreover $X_0$ is connected because so are the $\tau^n(\mathcal{R})$'s and because, given any two $n,m \in \mathbb{Z}$, 
the brick $B_{-i}^-$ is contained in $\tau^m(\mathcal{R}) \cap \tau^n(\mathcal{R})$ for a large enough $i \in \N$.
We need now Lemma 2.4 of \cite{Bonino:2009} apart from the fact that $F=\mathrm{Fix}(H)$ can meet $\mathrm{Bd}^+(\At)$ in the present paper. We bypass this minor difficulty as follows. Consider the strip $\widehat{\mathbb{A}} = \R \times [-1,2]$ and a brick decomposition $\widehat{\mathcal{D}}$ of $\widehat{\mathbb{A}} \setminus F$ extending $\widetilde{\mathcal{D}}_H$, that means such that any brick of $\widetilde{\mathcal{D}}_H$ is also a brick of $\widehat{\mathcal{D}}$. Recall that any set $X \subset S$ which is a union of bricks of $\widetilde{\mathcal{D}}_H$ is 
closed in $S$ hence $\overline{X} \subset X \cup F$; moreover $F$ only has isolated points 
so $\partial_{\widehat{\A}} \mathrm{Cl}_{\widehat{\A}} (X) = \partial \overline{X}$. We now let $X$ to be the 
union of $X_0$ with all the bounded connected components 
of $\widehat{\A} \setminus (F \cup X_0)$. Equivalently, $X$ is the union of $X_0$ with the connected components 
of $S \setminus X_0$ which are bounded and disjoint from $\mathrm{Bd}^+(\At)$. Thus $X$ is a connected  
union of bricks of $\widetilde{\mathcal{D}}_H \subset \widehat{\mathcal{D}}$ and \cite{Bonino:2009}[Lemma 2.4] 
tell us that $\partial_{\widehat{\A}} \mathrm{Cl}_{\widehat{\A}} (X) = \partial \overline{X}$ 
is a 1-dimensional submanifold of $\widehat{\A}$. Moreover $\mathrm{Bd}^-(\At) \subset \mathrm{Int}_{\widehat{\A}}(X)$ and $(\mathbb{R} \times \{2\}) \cap \mathrm{Cl}_{\widehat{\A}} (X) = \emptyset$ so there is a connected 
component $\tilde{J}$ of $\partial \overline {X}$ which is a line properly embedded 
in $\At$ and separating the boundary components of $\widehat{\A}$.
We have $\tau(X)=X$ hence $\tau(\partial \overline{X}) = \partial \overline{X}$ and then $\tau(\tilde{J}) = \tilde{J}$ 
since otherwise $\tilde{J}$ and $\tau(\tilde{J}) \subset \partial\overline{X}$ 
would be two disjoint properly embedded lines in $\At$ joining the two ends of $\At$, which is not possible since $\overline{X}$ is connected and contains $\mathrm{Bd}^-(\At)$ in its interior. It follows that $J=\Pi(\tilde{J})$ is an essential Jordan curve in $\A \setminus \mathrm{Bd}^-(\A)$ and we consider the subannulus $\B \subset \A$ 
bounded by $\mathrm{Bd}^-(\A) \cup J$. In order to show that $\B$ satisfies the required properties, it is 
enough to check that $H^{-1}(\tilde{J} \cap S)$ is included in the connected component of $\At \setminus \tilde{J}$ 
containing $\mathrm{Bd}^-(\At)$. Since $X$ is closed in $S$ and $S$ is an open subset of $\At$ one gets 
$\partial \overline{X} \subset \partial_S X \cup F$. Moreover $\partial_S X \subset \partial_S X_0$
and the closedness of $X_0$ in $S$ implies $\partial_S X_0 \subset \bigcup_{n \in \mathbb{Z}} 
\partial_S(\tau^n(\mathcal{R}))$. Hence any point $\tilde{z} \in \tilde{J} \cap S$ belongs to
$\partial_S (\tau^n(\mathcal{R})) = \tau^n(\partial_S \mathcal{R})$ for some $n \in \mathbb{Z}$ 
and consequently $$ H^{-1}(\tilde{z}) \in H^{-1}(\tau^n (\partial_S \mathcal{R})) 
= \tau^n (H^{-1}(\partial_S \mathcal{R})) \subset \tau^n(\mathrm{Int}(\mathcal{R})) 
\subset \mathrm{Int}(X).$$ We are done since $\mathrm{Int}(X)$ is connected, disjoint from $\tilde{J}$ 
and contains $\mathrm{Bd}^-(\At)$. \cqfd 

\ \\
Proposition \ref{p3} is a consequence of Lemmas \ref{l3}-\ref{l4} below.
\begin{lemme} \label{l3}
If $\mathcal{A}$ is bounded on the left and meets $\mathrm{Bd}^+(\At)$ then at least one of the two following 
assertions is true.
\begin{enumerate}
\item There exists an essential Jordan curve $J \subset \A \setminus \mathrm{Bd}(\A)$ such that $J \cap h(J) = \emptyset$. 
In particular the alternative (1') of Theorem \ref{t2} occurs.
\item There exists an arc $\tilde{\beta}$ crossing $\At$ such that, writing $W_r$ for the domain on the 
right of $\tilde{\beta}$, we have
\begin{itemize}
\item $\tau(W_r) \subset H(W_r) \subset W_r$; equivalently $G(W_r) \subset W_r \subset H^{-1}(W_r)$ since $G=\tau \circ H^{-1}$.
\item $H(\tilde{\beta}) \cap \tilde{\beta} = \tilde{\beta} \cap F$.
\end{itemize}
\end{enumerate}
\end{lemme}
\emph{Proof of Lemma \ref{l3}.} Suppose that $\mathcal{A}$ is bounded on the left 
and meets $\mathrm{Bd}^+(\At)$. As for Proposition \ref{p2}, we consider the strip 
$\widehat{\A} = \R \times [-1,2]$ and we endow the surface $\widehat{\A} \setminus F$  
with a brick decomposition $\widehat{\mathcal{D}}$ extending $\widetilde{\mathcal{D}}_H$. We let $X$ to be the union of $\mathcal{A}$ with the bounded connected components of $\widehat{\A} \setminus (F \cup \mathcal{A})$. Thus $X$ is a connected union of bricks of $\widetilde{\mathcal{D}}_H 
\subset \widehat{\mathcal{D}}$, it is bounded on the left and meets the two boundary components of $\At$. We have $\partial_{\widehat{\A}} \mathrm{Cl}_{\widehat{\A}} (X) = \partial \overline{X}$ and, according 
to \cite{Bonino:2009}[Lemma 2.4], this set is a 1-dimensional submanifold of $\widehat{\A}$.  
Hence there is a connected component $\Delta$ of $\partial \overline{X}$ which is a half-line properly embedded in $\At$, originating on $\mathrm{Bd}^-(\At)$ and meeting $\mathrm{Bd}^+(\At)$. We define $\tilde{\alpha}$ to 
be the arc obtained by following $\Delta$ from its origin to the first point where it intersects $\mathrm{Bd}^+(\At)$. Let $U_r$ be the domain on the right of $\tilde{\alpha}$ and observe that 
$\mathrm{Int}(X) \subset U_r$ because $\mathrm{Int}(X)$ is unbounded on the right and disjoint from $\tilde{\alpha} \subset \partial \overline{X}$. Let us show that $H(U_r) \subset U_r$. 
Remember that $\partial U_r = \tilde{\alpha} \subset \partial \overline  X \subset \partial_S X \cup F \subset \mathcal{A} \cup F$ so $$\partial H(U_r) = H(\tilde{\alpha}) \subset H(\tilde{\alpha} \cap S) \cup (\tilde{\alpha} \cap F)$$ 
and 
$$H(\tilde{\alpha} \cap S) \subset H(\mathcal{A}) \subset \mathrm{Int}(\mathcal{A}) \subset \mathrm{Int}(X) \subset U_r.$$ 
This gives $\partial H(U_r) \subset \overline{U_r}$ and consequently $H(U_r) \subset U_r$ because 
$H$ fixes the two ends of $\At$. 

Let us define  $V = \bigcup_{n \in \N} G^n(U_r)$. Recall that $G$ is a lift of $h$ or $h^{-1}$ without fixed point 
so Winkelnkemper's version of the Poincar\'e-Birkhoff theorem (\cite{Winkelnkemper:1988}) tell us that either $G$ is conjugate to $\tau$ or there exists an essential Jordan curve $J \subset \A$ such that $J \cap h^{\pm 1}(J) = \emptyset$. 
In the latter case, $J$ is certainly disjoint from $\mathrm{Bd}(\A)$ and we are done. We suppose now that $G$ is 
conjugate to $\tau$. It follows that $(\theta, \pm 1) < G(\theta,\pm 1)$ on $\mathrm{Bd}^{\pm}(\At)$ for every 
point $(\theta,\pm 1) \in \mathrm{Bd}(\At)$ since this is already known to be true for $(\theta,-1) \in \mathrm{Bd}^-(\At)$ hence we have $V = \bigcup_{0 \leq n \leq m} G^n(U_r)$ for some $m \in \N$. This implies $\overline{V} = \bigcup_{0 \leq n \leq m} G^n(\overline{U_r})$ so $\At \setminus \overline{V}$ has a single unbounded (on the left) connected component which we call $W_l$. It is a classical result of Ker\'ekj\'art\'o that any connected component $W$ of the intersection of two Jordan domains $U_1,U_2 \subset \R^2$ is also a Jordan domain, with frontier 
$\partial_{\R^2} W \subset \partial_{\R^2} U_1 \cup \partial_{\R^2} U_2$ (see \cite{Kerekjarto:1923} 
or \cite{LeCalvez/Yoccoz:1997}[Section 1]). One deduces that $\partial W_l$ is an arc 
crossing $\At$ and contained in $\bigcup_{0 \leq n \leq m} G^n(\tilde{\alpha})$. We let $\tilde{\beta} = \partial W_l$ 
and we write $W_r$ for the domain on the right of $\tilde{\beta}$. It follows from $G(V) \subset V$ that 
$W_l \subset G(W_l)$, i.e. $H(W_l) \subset \tau(W_l)$, and then $\tau(W_r) \subset H(W_r)$. Furthermore 
$H \circ G = G \circ H$ and $H(U_r) \subset U_r$ give $H(V) \subset V$ so $W_l \subset H(W_l)$ and then $H(W_r) \subset W_r$. It remains to check that $\tilde{\beta}$ and $H(\tilde{\beta})$ meet only in $F$. Just observe that $G(S)=S$ hence 
$$ \tilde{\beta} \cap S \subset \big ( \bigcup_{n \in \N} G^n(\tilde{\alpha}) \big) \cap S \subset \bigcup_{n \in \N} G^n(\partial_S \mathcal{A})$$ and consequently
$$H(\tilde{\beta} \cap S) \subset \bigcup_{n \in \N} G^n(H(\partial_S \mathcal{A})) \subset \bigcup_{n \in \N} G^n
(\mathrm{Int}(\mathcal{A})) \subset \bigcup_{n \in \N} G^n(U_r) = V \subset W_r.$$
\cqfd

All the arguments for proving the next lemma are already present in \cite{Bonino:2009}. Indeed one gets (1) 
below as \cite{Bonino:2009}[Lemma 3.6] provided the map $H^2$ appearing in 
\cite{Bonino:2009} is replaced with the map $H$ of the present paper. Item (2) is essentially obtained in the same way 
from \cite{Bonino:2009}[Lemma 3.7] and close ideas were previously used in \cite{Guillou:1994}. For convenience, we write a proof with some details omitted and refer to \cite{Bonino:2009} for more complete arguments.   
  
\begin{lemme} \label{l4}
Let $\tilde{\beta}$ and $W_r$ be as in Lemma \ref{l3}. 
\begin{enumerate}
\item If $\tau(\overline{W_r}) \not \subset W_r$, i.e. if $\tau(\tilde{\beta}) \cap \tilde{\beta} \neq \emptyset$, 
then the alternative (1') of Theorem \ref{t2} occurs.
\item If $\tau(\overline{W_r}) \subset W_r$, i.e. if $\tau(\tilde{\beta}) \cap \tilde{\beta} = \emptyset$, then 
the alternative (2') of Theorem \ref{t2} occurs.
\end{enumerate}
\end{lemme}
\emph{Proof of Lemma \ref{l4}.} (1) Suppose that $\tau(\tilde{\beta}) \cap \tilde{\beta} \neq \emptyset$.
In the following, any arc crossing $\At$ is oriented from its endpoint on $\mathrm{Bd}^-(\At)$ to its endpoint on $\mathrm{Bd}^+(\At)$. We write $\tilde{z}_0, \tilde{z}_1$ for respectively the endpoint of $\tilde{\beta}$ on 
$\mathrm{Bd}^-(\At)$ and on $\mathrm{Bd}^+(\At)$; we also let $\tilde{z}$ to be the first point of $\tilde{\beta}$ such that $\tau(\tilde{z}) \in \tilde{\beta}$. We first assume that $\tilde{z}$ and $\tau(\tilde{z})$ are met in this 
order on $\tilde{\beta}$. Since $\tau(\tilde{\beta})$ approaches $\tilde{\beta}$ from only one side, any two points in $\tilde{\beta} \cap \tau(\tilde{\beta})$ are met in the same order on $\tilde{\beta}$ and on $\tau(\tilde{\beta})$. Consequently we have $\{\tau(\tilde{z})\} = [\tilde{z_0},\tau(\tilde{z})]_{\tilde{\beta}} \cap \tau(\tilde{\beta})$ and the set $\tilde{\gamma} = [\tilde{z}_0, \tau(\tilde{z})]_{\tilde{\beta}} \cup \tau([\tilde{z}_0, \tilde{z}]_{\tilde{\beta}})$ is the frontier $\partial \Omega$ of some connected component 
$\Omega$ of $W_r \cap \tau(W_l)$, where $W_l$ denotes again the domain on the left of $\tilde{\beta}$. 
Let us show that $\B= \Pi(\overline{\Omega})$ is a subannulus of $\A$ as required. 
We let $\tilde{J} = [\tilde{z},\tau(\tilde{z})]_{\tilde{\beta}} \subset \At \setminus \mathrm{Bd}(\At)$. 
Clearly $\tilde{J}$ projects onto an essential loop $J=\Pi(\tilde{J}) \subset \A \setminus \mathrm{Bd}(\A)$
hence, in order prove that $\B$ is an essential subannulus of $\A$, one just need 
to show that the covering map $\Pi$ is one-to-one when restricted to $\overline{\Omega} \setminus [\tilde{z_0},\tilde{z}]_{\tilde{\beta}}$. Equivalently, one has to show 
$$\forall n \in \Z \setminus \{0\} \quad 
\tau^n(\overline{\Omega} \setminus [\tilde{z_0},\tilde{z}]_{\tilde{\beta}}) \cap 
(\overline{\Omega} \setminus [\tilde{z_0},\tilde{z}]_{\tilde{\beta}}) = \emptyset.$$
This is true for $n=\pm 1$ because $\Omega \subset W_r \setminus \tau(\overline{W_r})$ so 
$$\tau(\overline{\Omega}) \cap \overline{\Omega} \subset \tau(\overline{W_r}) \cap \partial \Omega 
= \tau([\tilde{z}_0,\tilde{z}]_{\tilde{\beta}}).$$
Moreover $\overline{\Omega}$ is a topological closed disc and it is then classical that the inclusion  
$\tau(\overline{\Omega}) \cap \overline{\Omega} \subset \partial \Omega$ implies 
$\tau^n(\overline{\Omega}) \cap \overline{\Omega} = \emptyset$ for every integer $n \not \in \{0,\pm 1\}$ 
(see for example the footnote in \cite{Bonino:2009}[p.1916] for details). 

It remains to be checked that $\B$ contains its image by $h$ or $h^{-1}$ and that $h(J) \cap J = J \cap N_0$.  
It is enough to show that $H(\tilde{J}) \subset \overline{\Omega}$ and that 
$H(\tilde{J}) \cap \tilde{J} = \tilde{J} \cap F$. For this last equality just recall 
that $H(\tilde{J}) \cap \tilde{J} \subset H(\tilde{\beta}) \cap \tilde{\beta} = \tilde{\beta} \cap F$.
Observe now that $\tilde{\beta} \cap \tau(\tilde{\beta}) \subset  \tilde{\beta} \cap H(\tilde{\beta}) \subset F$ so $\tau^k(\tilde{z}) \in F$ for any $k \in \Z$. We have $H(\tilde{J}) \cap \Omega \neq \emptyset$ 
since the H-image of a point of $\tilde{J} \cap S$ close to $\tilde{z}$ is a point of $W_r \cap \tau(W_l)$ 
close to $H(\tilde{z}) = \tilde{z}$ and since $\Omega$ is the only connected component of $W_r \cap \tau(W_l)$ 
having $\tilde{z}$ in its frontier. Finally $H(\tilde{J})$ lies entirely in $\overline{\Omega}$ 
since otherwise $H(\tilde{J} \setminus \{\tilde{z},\tau(\tilde{z})\}) \subset \overline{W_r} \cap 
\tau(\overline{W_l})$ would contain the point $\tau(\tilde{z}) = H(\tau(\tilde{z}))$, a contradiction. This proves the 
first item when $\tilde{z}$ precedes $\tau(\tilde{z})$ on $\tilde{\beta}$. In the other case, observe that we also 
have $\tau^{-1}(W_l) \subset H^{-1}(W_l) \subset W_l$ and replace in the above reasoning 
$\tilde{z},\tau(\tilde{z}),\tau,H,W_l,W_r$ with respectively $\tau(\tilde{z}),\tilde{z},\tau^{-1}, H^{-1}, W_r,W_l$.  

\ \\
(2) Let $X = \bigcap_{i \in \N} G^i(\tilde{\beta})$ and $X_n= \bigcap_{i=0}^n G^i(\tilde{\beta})$ ($n \geq 1$). Remark 
that $X = \emptyset$, i.e. $X_N= \emptyset$ for a least integer $N \geq 1$. Otherwise $G^{-1}$ induces a map from $X$ to $G^{-1}(X) \subset X$ and this map preserves the order on $X$ naturally provided by an orientation of 
$\tilde{\beta}$ because $G(\tilde{\beta})$ approaches $\tilde{\beta}$ from only one side. One would deduce 
that $G$ has a fixed point in $X$, a contradiction. 

If $N=1$, i.e. if $\tilde{\beta} \cap G(\tilde{\beta}) = \emptyset$, then $H(\tilde{\beta}) \cap 
\tau(\tilde{\beta}) = \emptyset$ and we simply define $\alpha = \Pi(\tilde{\beta})$. 
If $N \geq 2$ a suitable modification of $\tilde{\beta}$ allows to bring down the integer $N$ and then to 
reduce inductively to the easy case $N=1$. Let us give a few details. One has 
\begin{itemize}
\item[(i)] $\tau(X_{N-1}) \subset H(W_r) \subset W_r$, i.e. $G(X_{N-1}) \subset W_r \subset H^{-1}(W_r)$, 
\item[(ii)] $H(X_{N-1}) \subset W_r$,  
\item[(iii)] $X_{N-1} \cap \mathrm{Bd}(\At) = \emptyset$.
\end{itemize}
Property (i) follows from $G(X_{N-1}) \cap \tilde{\beta} = X_N= \emptyset$ and (ii) from 
$H(X_{N-1}) \cap \tilde{\beta} \subset H(G(\tilde{\beta})) \cap \tilde{\beta} 
= \tau(\tilde{\beta}) \cap \tilde{\beta} = \emptyset$. For (iii), just recall again that $G$ is fixed point free.    
In particular (i)-(ii) imply 
\begin{itemize} 
\item[(iv)] 
$\tau(X_{N-1}) \cap X_{N-1} = \emptyset = H(X_{N-1}) \cap X_{N-1}$. 
\end{itemize}
One can find an open neighbourhood $U \subset \At \setminus \mathrm{Bd}(\At)$ of the compact 
set $X_{N-1}$ in $\At$ which is so small that (i)-(iv) remain true when one replaces $X_{N-1}$ with $U$. Consider a finite covering $X_{N-1} \subset \bigcup_{i=1}^n \tilde{\alpha}_i$ where the $\tilde{\alpha}_i$'s are some connected components of $U \cap \tilde{\beta}$. For each $i=1,\cdots,n$, pick an arc $\tilde{\gamma}_i$ as follows: 
\begin{itemize} 
\item $\tilde{\gamma}_i$ lies entirely in $W_l \cap U$ except for its two endpoints points $\tilde{a}_i,\tilde{b}_i$ in $\tilde{\alpha}_i$; 
\item the arc $\tilde{\beta}_i \subset \tilde{\alpha}_i$ with the same endpoints $\tilde{a}_i,\tilde{b}_i$ as $\tilde{\gamma}_i$ is long enough to have 
$X_{N-1} \cap \tilde{\alpha}_i = X_{N-1} \cap (\tilde{\beta}_i \setminus \{\tilde{a}_i,\tilde{b}_i\})$.
\end{itemize}
Moreover these $\tilde{\gamma}_i$'s can be chosen pairwise disjoint hence we get a new arc $\tilde{\gamma}$ crossing $\At$ by removing from $\tilde{\beta}$ all the $\tilde{\beta}_i$'s and replacing them with the $\tilde{\gamma}_i$'s. 
Observe that $\tilde{\gamma} \cap X_{N-1} = \emptyset$ by the construction and, letting $W'_r$ be the domain on 
the right of $\tilde{\gamma}$, that $W_r \subset W'_r$. 
It is not difficult to prove that 
\begin{itemize}
\item[(a)] $\tau(W'_r) \subset H(W'_r) \subset W'_r$, 
\item[(b)] $\tau(\overline{W'_r}) \subset W'_r$,   
\item[(c)] $H(\tilde{\gamma}) \cap \tilde{\gamma} = \tilde{\gamma} \cap F$.
\end{itemize}
Moreover one checks that $\tilde{\gamma} \cap G(\tilde{\gamma}) \subset \tilde{\beta} \cap G(\tilde{\beta})$ 
hence $\bigcap_{i=0}^n G^i(\tilde{\gamma}) \subset X_n$ for every $n \geq 1$. Taking $n=N-1$ we obtain  
\begin{itemize}
\item[(d)] $\bigcap_{i=0}^{N-1} G^{i}(\tilde{\gamma}) \subset \tilde{\gamma} \cap X_{N-1} = \emptyset.$
\end{itemize} 
\, \cqfd

\section{Relationship with previous works} \label{section5}
Let us write $p_1$ for the projection $\At \to \R$, $(\theta,r) \mapsto \theta$. If $h$ 
is a homeomorphism of $\A$ isotopic to the identity and $H$ a lift of $h$ to $\At$, the 
\emph{horizontal displacement function} $D_H:\A \to \R$ is the continuous map defined 
by $D_H(z) = p_1 (H(\tilde{z})) - p_1(\tilde{z})$ where $z \in \A$ and $\tilde{z}$ is any point in $\Pi^{-1}(\{z\})$.

\subsection{Link with Franks' twist assumption} \label{A1}
The next result is Theorem 3.3 of \cite{Franks:1988a} (see also \cite{Franks:2006}).
\begin{thme} \label{tF}
Let $h:\A \to \A$ be a homeomorphism isotopic to the identity and with every point non wandering, and let 
$H: \At \to \At$ be a lift of $h$. If there exist points $\tilde{x},\tilde{y} \in \At$ such that 
$$\liminf_{n \to + \infty} \dfrac{p_1 (H^n(\tilde{y})) - p_1(\tilde{y})}{n} < 0 <  
\limsup_{n \to + \infty} \dfrac{p_1 (H^n(\tilde{x})) - p_1(\tilde{x})}{n}$$ 
then $h$ has at least two fixed points. 
\end{thme}
This roughly means that, for a conservative homeomorphism $h$ of $\A$, the presence of points turning 
clockwise and counterclockwise with linear speed under forward iteration guarantees the existence of two 
fixed points. The following consequence of Theorems \ref{t1}-\ref{t2} drops the 
assumption about the speed the points turn in $\A$ (see also Question \ref{q1} below). It also relax the generalized conservative assumption in Theorem \ref{tF} since if $\B \subset \A$ is a subannulus such that $h(\B) \varsubsetneqq \B$ then $\mathrm{Int}_{\A} (\B) \setminus h(\B)$ is a nonempty wandering open set for $h$ (similarly if $h^{-1}(\B) \varsubsetneqq \B$).   
\begin{thme} \label{t3}
Let $h:\A \to \A$ be a homeomorphism isotopic to the identity such that there is no essential subannulus 
$\B \subset \A$ containing $h(\B)$ or $h^{-1}(\B)$ as a proper subset, and let $H: \At \to \At$ be a lift of $h$. 
If there exist points $\tilde{x},\tilde{y} \in \At$ such that 
$$ \liminf_{n \to + \infty} p_1 (H^n(\tilde{y})) = - \infty \quad \text{ and }  \quad \limsup_{n \to +\infty} p_1 (H^n(\tilde{x})) = + \infty $$ 
then $h$ has at least two fixed points; more precisely the Nielsen class $\Pi(\mathrm{Fix}(H))$ contains at least 
two points. 
\end{thme}
\emph{Proof of Theorem \ref{t3}.} Let us check separately the first assertion. Suppose that $h$ has at most one fixed point. Consider an arc $\alpha$ given by Theorem \ref{t1} and a connected component $\tilde{\alpha}$ of $\Pi^{-1}(\alpha)$. Thus $\tilde{\alpha}$ is an arc crossing $\At$ and it has no point of transverse intersection with $H(\tilde{\alpha})$ because the same is true for $\alpha$ and $h(\alpha)$. It follows that one of the two connected components of $\At \setminus \tilde{\alpha}$, call it $W$, contains its image by $H$. Consequently either every $\tilde{z} \in W$ has its forward $H$-orbit bounded on the right or every $\tilde{z} \in W$ has its forward $H$-orbit bounded on the left. Since any $\tilde{z} \in \At$ has a translate $\tau^n(\tilde{z})$ in $W$ and since 
$H \circ \tau = \tau \circ H$, the last sentence remains true if one replaces ``$\tilde{z} \in W$"  with ``$\tilde{z} \in \At$". The first part of the theorem is proved. Suppose now that the Nielsen class 
$N_H = \Pi(\mathrm{Fix}(H))$ contains at most one point. Given a positive integer $n$, we define 
$\Ab = \R/n\Z \times [-1,1]$ and we consider the two natural covering maps $\breve{\Pi}: \At \to \Ab$ 
and $q:\Ab \to \A$. Of course $\Ab$ is homeomorphic to $\A$ with boundary components $\mathrm{Bd}^{\pm}(\Ab) = 
\breve{\Pi}(\mathrm{Bd}^{\pm}(\At)) = q^{-1}(\mathrm{Bd}^{\pm}(\A))$ and we write $\breve{h}, \breve{\tau}$ 
for the homeomorphisms of $\Ab$ induced by respectively $H,\tau$. Thus $\breve{h}$ is lifted by $H$ and 
$\breve{h}$ is a lift of $h$, that means $\breve{h} \circ \breve{\Pi} = \breve{\Pi} \circ H$ and $h \circ q = q \circ \breve{h}$. Moreover one can choose $n$ so large that $\breve{\Pi}(\mathrm{Fix}(H)) = \mathrm{Fix}(\breve{h})$ because the map $\tilde{z} \mapsto p_1(H(\tilde{z})) - p_1(\tilde{z})$ ($\tilde{z} \in \At$) is bounded. 
It is well known that if a Nielsen class of $h$ has finitely many fixed points then the sum of their Lefschetz indexes 
equals $\chi(\A)=0$ so $N_H$ is either empty or contains a single point $z_0$ with Lefschetz index $\mathrm{Ind}_h(z_0)=0$ 
\footnote{A direct argument is the following. Pick $\tilde{z}_0 \in \Pi^{-1}(\{z_0\})$. We 
first assume that $z_0 \not \in \mathrm{Bd}(\A)$ and we choose a rectangle $R= [a,a+1] \times[-1,1]$ containing $\tilde{z}_0$ in its interior. We have then $R \cap \mathrm{Fix}(H) = \{\tilde{z}_0\}$ and it is easily seen that the 
index of the Jordan curve $J= \partial_{\R^2} R$ w.r.t $H$ is equal to $0$. It follows that 
$0=\mathrm{Ind}_H(\tilde{z}_0) = \mathrm{Ind}_h(z_0)$. If $z_0$ lies on the boundary of $\A$, say $z_0 \in \mathrm{Bd}^+(\A)$, extend $H$ to a homeomorphism $H_{\ast}$ of the ``double strip" $\R \times [-1,3]$ by letting $H_{\ast} = S \circ H \circ S$ on $\R \times [1,3]$ where $S(\theta,r)=(\theta,2-r)$. One gets as above $0 = \mathrm{Ind}_{H_{\ast}}(\tilde{z}_0) = 2 \mathrm{Ind}_h(z_0)$.}. Consequently $\mathrm{Fix}(\breve{h})$ consists in a single Nielsen class 
of $\breve{h}$ which is either empty or contains exactly $n$ points; in the latter case, these $n$ fixed points 
have index $0$ w.r.t. $\breve{h}$ and they all lie either in the same boundary component of $\Ab$ or in 
$\Ab \setminus \mathrm{Bd}(\Ab)$. Hence we can apply Theorem \ref{t2} to $\breve{h}: \Ab \to \Ab$. If (2') occurs we 
get in particular an arc $\breve{\alpha}$ crossing $\Ab$ which has no point of transverse intersection with $\breve{h}(\breve{\alpha})$ and we conclude as above by lifting $\breve{\alpha}$ to $\At$. We end by showing that the existence 
of a subannulus $\Bb \subset \Ab$ given by (1') contradicts our hypotheses. We assume that $\breve{h}(\Bb) \varsubsetneqq \Bb$ and that the boundary of $\Bb$ is the union of $\mathrm{Bd}^-(\Ab)$ together with a Jordan 
curve $\breve{J} \subset \Ab \setminus \mathrm{Bd}^-(\Ab)$, the other cases being similar. 
We consider the annulus $\Ab' = \R/n\Z \times [-1,2]$ and we let $U$ to be the connected component of 
$\Ab' \setminus \bigcup_{0 \leq i \leq n-1} \breve{\tau}^i(\Bb)$ containing the upper boundary $\R/n\Z \times \{2\}$. 
According to Ker\'ekj\'art\'o's result mentioned in the proof of Lemma \ref{l3}, the set $\partial_{\Ab'} U$ 
is a Jordan curve contained in $\bigcup_{0 \leq i \leq n-1} \breve{\tau}^i(\breve{J}) \subset \Ab \setminus 
\mathrm{Bd}^-(\Ab)$. For any $\breve{z} \in \breve{J} \setminus \mathrm{Fix}(\breve{h})$ and any $i=0, \cdots, n-1$ one 
has $\breve{h} \circ \breve{\tau}^i (\breve{z}) = \breve{\tau}^i \circ \breve{h} (\breve{z}) \in \breve{\tau}^i 
(\mathrm{Int}_{\Ab}(\Bb))$. Moreover $\bigcup_{i=0}^{n-1} \breve{\tau}^i (\mathrm{Int}_{\Ab}(\Bb))$ is connected, 
disjoint from $\partial_{\Ab'} U$ and contains $\mathrm{Bd}^-(\Ab)$ hence $\breve{h}(\partial_{\Ab'} U \setminus 
\mathrm{Fix}(\breve{h})) \subset \Ab' \setminus \mathrm{Cl}_{\Ab'} (U) \subset \Ab$. Finally $\breve{\tau} (\partial_{\Ab'} U) = \partial_{\Ab'} U$ so $J = q(\partial_{\Ab'} U)$ is an essential Jordan curve in $\A \setminus \mathrm{Bd}^-(\A)$. 
Letting $\B$ to be the subannulus of $\A$ bounded by $J \cup \mathrm{Bd}^-(\A)$, one deduces that 
$h(\B) \varsubsetneqq \B$, a contradiction. \cqfd 
\begin{rem}
\begin{itemize}
\item The same argument could be used with \cite{Guillou:1994}[Th\'eor\`eme 5.1] to find one fixed point 
of $h$ under the assumptions of Theorem \ref{t3}. Another result by Franks (\cite{Franks:1988c}[Corollary 2.3]) with a somewhat different generalized conservative assumption ($h$ is assumed to be chain transitive) also ensures the existence of 
one fixed point for $H$ if there exist $\tilde{x}, \tilde{y}$ as in Theorem \ref{t3}. 
\item Observe that in the proof of Theorem \ref{t3} it was enough to find an arc $\tilde{\alpha}$ crossing $\At$ 
which has no point of transverse intersection with $H(\tilde{\alpha})$, even if it does not project onto an arc of $\A$. 
In contrast, the full properties of $\alpha$ are needed for Theorem \ref{t4} below.    
\end{itemize}
\end{rem}
\begin{question} \label{q1} 
I do not know an example of conservative homeomorphism $h:\A \to \A$ isotopic to the identity with a lift $H$ satisfying 
the generalized twist property of Theorem \ref{t3} and not the one in Theorem \ref{tF}, even if the word conservative is understood in a generalized sense. Such an example is easily constructed if one drops the conservative assumption. One can formulate this problem by using the rotation set of $H$. We recall this notion is a well-known adaptation of a similar concept introduced by M. Misiurewicz and K. Ziemian for maps 
of tori (\cite{Misiurewicz/Ziemian:1989}). 
Given a lift $H:\At \to \At$ of $h$, the rotation set $\rho(H) \subset \R$ may be defined by deciding that $r \in \rho(H)$ if there exist a sequence $(n_k)_{k \in \N}$ of positive integers and a sequence 
$(\tilde{z}_k)_{k \in \N}$ of points of $\At$ such that 
$$\lim_{k \to + \infty} n_k = + \infty \text{ and } r = \lim_{k \to + \infty} \dfrac{p_1(H^{n_k}(\tilde{z}_k)) - p_1(\tilde{z}_k)}{n_k}.$$
It is a compact interval, say $\rho(H) = [a,b]$ with possibly $a=b$. Equivalently, one has $r \in \rho(H)$ if $r = \int_{\A} D_H(z) d\mu(z)$ for some $h$-invariant Borel probability measure $\mu$ on $\A$. It follows that the endpoints of $\rho(H)$ are the rotation number of some points, that means
$$a = \lim_{n \to + \infty} \frac{p_1(H^{n}(\tilde{z})) - p_1(\tilde{z})}{n}, \quad 
b = \lim_{n \to + \infty} \frac{p_1(H^{n}(\tilde{z}')) - p_1(\tilde{z}')}{n} $$ for some 
$\tilde{z}, \tilde{z}' \in \At$.
Hence there exist $\tilde{x}, \tilde{y}$ as in Theorem \ref{tF} if and only if 
$\rho(H)$ contains $0$ in its interior. 
Moreover one checks that if $0 \not \in \rho(H)$ then all the forward orbits of $H$ go to the same end of $\At$ hence 
our initial question is essentially the same as the following.
Can a conservative homeomorphism $h$ of $\A$ be lifted by $H:\At \to \At$ having a forward orbit unbounded 
on the left and such that $\rho(H) = [0,b]$, $b \geq 0$ ?       
\end{question} 
\subsection{Link with Conley-Zehnder theorem in the annulus} \label{A2} 
We write $\lambda$ for the Lebesgue measure on $\A$ and $\tilde{\lambda}$ for the lift of $\lambda$ to $\At$. Recall that the measure $\tilde{\lambda}$ is invariant by $\tau$ 
and that $\tilde{\lambda}(X) = \lambda(\Pi(X))$ for any Borelian set $X$ contained in a fundamental domain of $\Pi$. If $h$ is supposed to preserve $\lambda$ then 
$\tilde{\lambda}$ is also invariant by any lift $H:\At \to \At$ of $h$. We first state a slight complement to Theorem \ref{t2}.

\begin{prop} \label{p4} 
Let $h:\A \to \A$ be a homeomorphism isotopic to the identity. We suppose that $h$ satisfies the assumptions (i)-(iii) in Theorem \ref{t2} (in particular this holds if $h$ has at most one fixed point) and furthermore that $h$ preserves the measure $\lambda$. Then the alternative (2') in Theorem \ref{t2} occurs and (1') does not; 
moreover the arc $\alpha$ in (2') can be chosen such that $\lambda(\alpha) = 0$. 
\end{prop}
\emph{Proof of Proposition \ref{p4}.} There is no subannulus $\B \subset \A$ as in Theorem \ref{t2} since otherwise $h$ would have a nonempty wandering open set, contradicting $\lambda(\A) < \infty$. Let us remark now that our constructions can be improved in order to get $\lambda(\alpha) = 0$ in Theorem \ref{t2}. First, the brick decomposition 
$\widetilde{\mathcal{D}}_H = \{B_i\}_{i \in I}$ given by Lemma \ref{l1} can be chosen in such a way $\tilde{\lambda}(\bigcup_{i \in I} \partial B_i) = 0$; it suffices for example to deal only with bricks decompositions of $S$ having polygonal bricks. Using now the notation from Lemma \ref{l3} and its proof, 
observe that the arc $\tilde{\beta}$ is a subset of $\bigcup_{n \in \N} G^n(\tilde{\alpha})$ where $\tilde{\alpha}$ is itself a subset of $(\bigcup_{i \in I} \partial B_i) \cup F$. Since $F$ is either empty or a countable set and since 
$G$ preserves the measure $\tilde{\lambda}$ one obtains $\tilde{\lambda}(\tilde{\beta})=0$. Keeping finally 
the notation from Lemma \ref{l4} (2) and its proof, it remains to see that the modification 
of $\tilde{\beta}$ into $\tilde{\gamma}$ (if needed, i.e. if $N \geq 2$) can be performed in such a way 
that $\tilde{\lambda}(\tilde{\gamma}) = 0$. It is enough to observe that the arcs $\tilde{\gamma}_i$, $i=1,\cdots,n$, can be chosen such that $\tilde{\lambda}(\tilde{\gamma}_i)=0$ because, for a given $i$, there are 
uncountably many pairwise disjoint arcs $\tilde{\gamma}_i$ with the required properties while only countably many of them 
can have positive $\tilde{\lambda}$-measure (as an alternative argument, one can get $\tilde{\lambda}(\tilde{\gamma}_i)=0$ by constructing $\tilde{\gamma}_i$ piecewise linear; this is possible by using that densely many points in $\tilde{\alpha}_i$ are accessible from $W_l$ by a straight segment). After $N-1$ such modifications one reduces to the case $\tilde{\beta} \cap G(\tilde{\beta})= \emptyset$ with $\tilde{\lambda}(\tilde{\beta})=0$, so 
$\alpha=\Pi(\tilde{\beta})$ is an arc as described in Theorem \ref{t2} with furthermore $\lambda(\alpha)=0$. \cqfd

\ \\
The next lemma gives an interesting interpretation for the mean horizontal displacement of $H$. This seems to be 
well known although I found it explicitly only in \cite{Beguin/Crovisier/LeRoux:2006}. 
\begin{lemme} \label{l5}
Let $h:\A \to \A$ be a homeomorphism isotopic to the identity, preserving the measure $\lambda$, and 
let $H: \At \to \At$ be a lift of $h$. Consider an arc $\alpha$ crossing $\A$ such that $\lambda(\alpha)=0$ and 
a connected component $\tilde{\alpha}$ of $\Pi^{-1}(\alpha)$. Then the mean horizontal displacement 
$\int_{\A} D_H(z) \, d\lambda(z)$ is equal to the algebraic area (w.r.t the measure $\tilde{\lambda}$) between $\tilde{\alpha}$ and its image $H(\tilde{\alpha})$.
\end{lemme}  
Under the assumptions of the above lemma, recall from Birkhoff's ergodic theorem that $\lambda$-almost every 
$z \in \A$ has a rotation number $\rho_H(z) \in \R$ defined by $\rho_H(z) = \lim_{n \to + \infty} 
\frac{1}{n} \big (p_1(H^n(\tilde{z})) - p_1(\tilde{z}) \big)$ for any $\tilde{z} \in \Pi^{-1}(\{z\})$ and that the integrable map $z \mapsto \rho_H(z)$ satisfies  
$\int_{\A} \rho_H(z) \, d\lambda(z) = \int_{\A} D_H(z) \, d\lambda(z)$. Hence Lemma \ref{l5} is given by Proposition 5.2 
of \cite{Beguin/Crovisier/LeRoux:2006} when the arc $\alpha$ is the projection of a vertical segment $\{\theta\} \times [-1,1] \subset \At$ (the boundedness assumption on $D_H$ appearing in \cite{Beguin/Crovisier/LeRoux:2006} is needed only when working in the open annulus and is always satisfied in our compact framework). To obtain exactly the above statement, use the Oxtoby-Ulam theorem on homeomorphic measures (\cite{Oxtoby/Ulam:1941}) to construct a homeomorphism $g:\A \to \A$ isotopic to the identity, preserving the measure $\lambda$ and mapping $\alpha$ onto a radial segment $\alpha' = \{e^{2i\pi\theta}\} \times [-1,1]$. Let  $G:\At \to \At$ be the lift of $g$ mapping $\tilde{\alpha}$ onto 
$\tilde{\alpha}' = \{\theta\} \times [-1,1]$ and define $H' = G \circ H \circ G^{-1}$ which is a lift of 
$h'= g \circ h \circ g^{-1}$. Of course the algebraic area between $\tilde{\alpha}$ and 
$H(\tilde{\alpha})$ is the same as the one between $G(\tilde{\alpha}) = \tilde{\alpha}'$ and 
$G(H(\tilde{\alpha})) = H'(\tilde{\alpha}')$. One concludes because  
$$\int_{\A} D_H(z) \, d\lambda(z) = \int_{\A} D_{H'}(z) \, d\lambda(z).$$

The next result is essentially Theorem 2 in Flucher's paper \cite{Flucher:1990} restricted to the case of the 
Lebesgue measure. 
\begin{thme} \label{t4}
Let $h: \A \to \A$ be a homeomorphism isotopic to the identity and preserving the measure $\lambda$. 
If $h$ admits a lift $H:\At \to \At$ with vanishing mean horizontal displacement, i.e. with $\int_{\A} D_H(z) \, d\lambda(z) = 0$, then $h$ possesses at least two fixed points; more precisely the Nielsen class $\Pi(\mathrm{Fix}(H))$ contains at least two points. 
\end{thme}
\emph{Proof of Theorem \ref{t4}.} Suppose first that $h$ has at most one fixed point. Consider an arc $\alpha$ given by 
Proposition \ref{p4} and a connected component $\tilde{\alpha}$ of $\Pi^{-1}(\alpha)$. For any lift $H$ of $h$, 
the arc $H(\tilde{\alpha})$ lies entirely in a connected component of $\At \setminus \tilde{\alpha}$ except 
possibly for one point in $\tilde{\alpha}$. The algebraic area between and $\tilde{\alpha}$ and 
$H(\tilde{\alpha})$ is then nonzero and so is $\int_{\A} D_H(z) \, d\lambda(z)$ by Lemma \ref{l5}. 
We suppose now that $H$ is a lift of $h$ such that $\Pi(\mathrm{Fix}(H))$ contains at most one fixed point. 
We consider a $n$-folded covering $\Ab$ of $\A$ and $\breve{h}:\Ab \to \Ab$ as in the proof of Theorem \ref{t3}. 
The homeomorphism $\breve{h}$ preserves the measure $\breve{\lambda}$ obtained by lifting $\lambda$ to the annulus $\Ab$. Observe that Proposition \ref{p4} applies with $\breve{h}, \breve{\lambda}$ instead of $h,\lambda$ (with the same proof since $H,\tilde{\lambda}$ are also the lifts of $\breve{h},\breve{\lambda}$) and consider an arc $\breve{\alpha} \subset \Ab$ obtained in this way. Defining the horizontal displacement function on $\Ab$ by $\breve{D}_H(\breve{z}) = 
\frac{1}{n} \big (p_1(H(\tilde{z})) - p_1(\tilde{z}) \big)$ (where $\breve{z} \in \Ab$ and $\tilde{z} \in \breve{\Pi}^{-1}(\{\breve{z}\})$) one checks that Lemma \ref{l5} is still valid with $\breve{h}, \breve{\alpha}, \breve{\lambda}, \breve{D}_H$ instead of $h,\alpha, \lambda, D_H$ and we get as above 
$$0 \neq \int_{\Ab} \breve{D}_H(\breve{z}) \,d\breve{\lambda}(\breve{z}) = \int_{\A} D_H(z) \,d\lambda(z).$$
\, \cqfd

\ \\
\textbf{Acknowledgements.} I would like to thank Sylvain Crovisier for several conversations and for bringing to my attention Proposition 5.2 of \cite{Beguin/Crovisier/LeRoux:2006}. I also thank the referee for her/his careful reading of the manuscript.

\bibliographystyle{plain}
\bibliography{/home2/bonino/biblio}
\end{document}